\documentclass[reqno,10pt]{amsart}%
\usepackage{color}
\usepackage{latexsym,amssymb,amsfonts}
\usepackage{mathrsfs}
\usepackage{amsmath}
\usepackage{graphicx}%

\renewcommand{\b}[1]{\mbox{\boldmath $#1$}}

\newtheorem{lemma}{\quad\bf Lemma}[section]

\newtheorem{corollary}{\quad\bf Corollary}
\newtheorem{remark}{\quad\bf Remark}

\newtheorem{theorem}{\quad\bf Theorem}[section]

\numberwithin{equation}{section} \theoremstyle{plain}
\theoremstyle{definition}

\theoremstyle{remark}
\newtheorem{example}{Example}[section]

\begin{document}
\author{Jak\v sa Cvitani\'c}
\address{Caltech, M/C 228-77, 1200 E. California Blvd. Pasadena, CA 91125, USA.}
\email{cvitanic@hss.caltech.edu}
\author{Robert Liptser}
\address{Department of Electrical Engineering-Systems, Tel Aviv University, 69978 Tel
Aviv, Israel} \email{liptser@eng.tau.ac.il}
\author{Boris Rozovskii }
\address{Department of Mathematics, University of Southern California, Los Angeles, CA
90089-1113, USA } \email{rozovski@math.usc.edu}
\thanks{The research of J. Cvitani\'c was supported in part by the National Science
Foundation, under Grant NSF-DMS-00-99549 and 04-03575.}
\thanks{The research of B.L. Rozovskii was supported in part by the Army Research
Office and the Office of Naval Research under the grants
DAAD19-02-1-0374 and N0014-03-0027.}

\date{}
\title[Filtering random-time observations]
{A filtering approach to tracking volatility from prices observed at random times}

\maketitle

\begin{abstract}
This paper is concerned with nonlinear filtering of the
coefficients in asset price models with stochastic volatility.
More specifically, we assume that the asset price process
$
S=(S_{t})_{t\geq0}
$
is given by
\[
dS_{t}=r(\theta_{t})S_{t}dt+v(\theta_{t})S_{t}dB_{t},
\]
where $B=(B_{t})_{t\geq0}$ is a Brownian motion, $v$ is a positive
function, and $\theta=(\theta_{t})_{t\geq0}$ is a c\'{a}dl\'{a}g
strong Markov process. The random process $\theta$ is
unobservable. We assume also that the asset price $S_{t}$ is
observed only at random times $0<\tau_{1}<\tau_{2}<\ldots.$ This
is an appropriate assumption when modelling high frequency
financial data (e.g., tick-by-tick stock prices).

In the above setting the problem of estimation of $\theta$ can be
approached as a special nonlinear filtering problem with
measurements generated by a multivariate point process
$(\tau_{k},\log S_{\tau_{k}})$. While quite natural, this problem
does not fit into the standard diffusion or simple point process
filtering frameworks and requires more technical tools. We derive
a closed form optimal recursive Bayesian filter for $\theta_{t}$ ,
based on the observations of $(\tau_{k},\log
S_{\tau_{k}})_{k\geq1}$. It turns out that the filter is given by
a recursive system that involves only deterministic
Kolmogorov-type equations, which should make the numerical
implementation relatively easy.
\end{abstract}

\section{\bf Introduction}
\label{Sec-1}

In the classical Black-Scholes model for financial markets, the
stock price $S_{t}$ is modelled as a Geometric Brownian motion,
that is, with diffusion coefficient equal to $\sigma S_{t}$, where
``volatility'' $\sigma$ is assumed to be constant. The volatility
parameter is the most important one when it comes to option
pricing; consequently, many researchers have generalized the
constant volatility model to so-called stochastic volatility
models, where $\sigma_{t}$ is itself random and time dependent.
There are two basic classes of models: complete and incomplete. In
complete models, the volatility is assumed to be a functional of
the stock price; in incomplete models, it is driven by some other
source of noise that is possibly correlated with the original
Brownian motion. In this paper we study a particular incomplete
model in which the volatility process is independent of the
driving Brownian motion process. This has the economic
interpretation of the volatility being influenced by market,
political, financial, and other factors that are independent of
the  ``systematic risk'' (the Brownian motion process) associated
with the particular stock price under study. It is also close in
spirit to the way traders think about volatility -- as a parameter
that changes with time and whose future value in a given period of
interest has to be estimated/predicted. They need an estimate of
the volatility to decide how they will trade in financial markets,
especially derivatives markets. In fact, the notion of volatility
is so important to traders that they even quote option prices in
volatility units rather than in dollars (or some other currency).
Investment banks also depend on modelling future volatility in
order to price custom-made financial products, whose payoff
depends on the future path of the underlying stock price.
Recently, new contracts have been developed that directly trade
the volatility itself (volatility swaps, for example). We plan to
address the issue of pricing options within the framework of our
model in future research.

Estimating volatility from observed stock prices is not a trivial
task in either complete or incomplete models, in part because the
prices are observed at discrete, possibly random time points.
Since volatility itself is not observed, it is natural to apply
filtering methods to estimate the volatility process from
historical stock price observations. Nevertheless, this has only
recently been investigated in continuous-time models, in
particular by Frey and Runggaldier \cite{FR}. See Runggaldier
\cite{R} for an up-to-date survey. See also Elliott et al
\cite{EHJ} for a discrete-time approach with equally spaced
observations, Gallant and Tauchen \cite{GT} for an approximating
algorithm in continuous time, Malliavin and Mancino \cite{MM} for
a nonparametric approach, as well as Fouque et al. \cite{FPS},
Rogers and Zane \cite{RZ}, and Kallianpur and Xiang \cite{KX} for
still other approaches. There is also a rich econometrics,
time-series literature on ARCH-GARCH models of stochastic
volatility, that presents an alternative way to model and estimate
volatility; see Gourieroux \cite{Gou} for a survey.

Our paper was prompted by Frey and Runggaldier \cite{FR}. Like
that paper, we assume that the asset price process
$S=(S_{t})_{t\geq0}$ is given by
\[
dS_{t}=r(\theta_{t})S_{t}dt+v(\theta_{t})S_{t}dB_{t},
\]
where $B=(B_{t})_{t\geq0}$ is a Brownian motion, $v$ is a positive
function, and $\theta=(\theta_{t})_{t\geq0}$ is a c\'{a}dl\'{a}g
strong Markov process. The ``volatility'' process $\theta$ is
unobservable, while the asset price $S_{t}$ is observed only at
random times $0<\tau_{1}<\tau_{2}<\ldots$ This assumption is
designed to reflect the discrete nature of high frequency
financial data such as tick-by-tick stock prices. The random time
moments $\tau_{k}$ can be interpreted as instances at which a
large trade occurs or at which a market maker updates his quotes
in reaction to new information  (see Frey \cite{F} ). Hence, it is
natural to assume that $(\tau_{k})_{k\geq1}$ might also be
correlated with $\theta.$

In the above setting the problem of volatility estimation can be
regarded as a special nonlinear filtering problem.

Frey and Runggaldier \cite{FR} derive a Kallianpur-Striebel type
formula (see e.g. \cite{KalStr}) for the optimal mean-square
filter for $\theta_{t}$ based on the observations of
$S_{\tau_{1}},S_{\tau_{2}},...$ for all $\tau_{k}\leq t$ and
investigate Markov Chain approximations for this formula. We
extend this result in that we derive the exact filtering equations
for $\theta_{t}$ that allow us to compute the conditional
distribution of $\theta_{t}$ given $S_{\tau_{1}\wedge t}$,
$S_{\tau_{2}\wedge t}$,\ldots. Moreover, our framework includes
general random times of observations, not just doubly stochastic
Poisson processes.

We remark that while being natural, the Frey and Runggaldier model
adopted in this paper does not quite fit into the standard
diffusion or simple point process filtering frameworks (cf.
\cite{LSII}, \cite{KrZa}, \cite{Roz1}) and requires more technical
tools. In particular, the general filtering theory for diffusion
processes requires that the diffusion coefficient of the
observation process does not depend on the state process, while in
our case the presence of $\theta_{t}$ in the diffusion coefficient
is crucial. The standard" filtering theory for point processes is
also not applicable in the present setting since the observation
process $(\tau_{i},S_{\tau_{i}})_{i\geq1}$ is a multivariate
process (see also Remark \ref{r1}).

It turns out that the resulting filtering equations are simpler
than their counterparts in the case of continuous observations. In
the latter case, the nonlinear filters are described by infinite
dimensional stochastic differential equations. For example, if
$\theta_{t}$ is a diffusion process, the filtering equations
(e.g., Kushner filter or Zakai filter) are given by stochastic
partial differential equations (see, e.g., \cite{Roz1}). In
contrast, in our setting, the filtering equation can be reduced to
a recursive system of linked \textit{deterministic} equations of
Kolmogorov's type. Therefore, the numerical implementation of the
filter is much simpler (see the follow up paper \cite{CRZ}).

We describe the model in Section 2, state the main results and
examples in Section 3, provide the proofs in Section 4, and
present more detailed examples in Section 5.

\section{\bf Mathematical model}
\label{Sec-2}

\subsection{Risky asset and observation times}

Let us fix a probability space $(\Omega,\mathcal{F},\mathsf{P})$
equipped with a filtration $\mathbf{F}=(\mathcal{F}_{t})_{t\geq0}$
that satisfies the usual conditions (see, e.g. \cite{LSMar}). All
random processes considered in the paper are assumed to be defined
on $(\Omega,\mathcal{F},\mathsf{P})$ and adapted to $\mathbf{F}$.

It is assumed that there is a risky asset with the price process
$S=(S_{t})_{t\geq0}$ given by the It\^{o} equation
\begin{equation}\label{2.1b}
dS_{t}=r(\theta_{t})S_{t}dt+v(\theta_{t})S_{t}dB_{t},
\end{equation}
where $B=(B_{t})_{t\geq0}$ is a standard Brownian motion and
$\theta =(\theta_{t})_{t\geq0}$ is a c\'{a}dl\'{a}g Markov
jump-diffusion process in $\mathbb{R}$ with the generator
$\mathcal{L}$. To simplify the discussion, it is assumed that
$r(x)$ and $v(x)$ are measurable bounded functions on
$\mathbb{R}$, the initial condition $S_{0}$ is constant, and
$v(x)$ and $S_{0}$ are positive.

The process $(\theta_{t})_{t\geq0}$ is referred to as the
\emph{volatility process}. It is unobservable, and the only
observable quantities are the values of the log-price process
$X_{t}=\log S_{t}$ taken at stopping times $(\tau_{k})_{k\geq0}$,
so that $\tau_{0}=0,\tau_{k}<\tau_{k+1}$ if $\tau _{k}<\infty,$
and $\tau_{k}$ $\uparrow\infty$ as $k\uparrow\infty.$

In accordance with \eqref{2.1b}, the log-price process is given by
\[
X_{t}=\int_{0}^{t}\Big(r(\theta_{s})-\frac{1}{2}v^{2}(\theta_{s}
)\Big)ds+\int_{0}^{t}v(\theta_{s})dB_{s}.
\]
For notational convenience, set $X_{k}:=X_{\tau_{k}}.$ Thus, the
observations are given by the sequence
$(\tau_{k},X_{k})_{k\geq0}$.

\begin{remark}\label{rem-1}
{\rm (Note on the reading sequence.) The reader interested
primarily in applying our results to real data can focus her
attention on Example \ref{examplechain}, which appears to be the
most practical model to work with. That example provides
self-contained formulas for estimating the conditional {\rm
(}filtering{\rm )} distribution of the volatility process. We
report on the numerical results related to this example in the
follow-up paper \cite{CRZ}.}
\end{remark}

Clearly, the observation process $(\tau_{k},X_{k})_{k\geq0}$ is a
multivariate (marked) point process (see, e.g. \cite{JS},
\cite{Last}) with the counting measure
\[
\mu(dt,dy)=\sum_{k\geq1}\mathbf{I}_{\{ \tau_{k}<\infty\} }
\delta_{{\{\tau_{k},X_{k}\}}}(t,y)dtdy,
\]
where $\delta_{{\{\tau_{k},X_{k}\}}}$ is the Dirac delta-function
on $\mathbb{R}_{+}\times\mathbb{R}$.

We introduce two filtrations related to
$(\tau_{k},X_{k})_{k\ge0}$: $(\mathcal{G}(n))_{n\ge0}$ and
$(\mathcal{G}_{t})_{t\ge0}$, where

- $\mathcal{G}(n):=\sigma\{(\tau_{k},X_{k})_{k\leq n}\}$,

- $\mathcal{G}_{t}:=\sigma(\mu([0,r]\times\Gamma):r\leq
s,\Gamma\in \mathcal{B}(\mathbb{R})),$ where
$\mathcal{B}(\mathbb{R})$ is the Borel $\sigma$-algebra on
$\mathbb{R.}$

\smallskip
\noindent
It is a standard fact (see III.3.31 in \cite{JS}) that
\begin{equation*}
\mathcal{G}_{\tau_{k}}=\mathcal{G}(k),\ k=0,1\ldots
\end{equation*}
and $\{\tau_{k}\}$ is a system of stopping times with respect to
$(\mathcal{G}_{t})_{t\geq0}$.

\begin{remark}\label{r1}
{\rm Although $\mathcal{G}_{\tau_{k}}$ contains all the relevant
information carried by the observations obtained up to time
$\tau_{k},$  the filtration $\big(\mathcal{G}_{t}\big)_{t\geq 0}$
provides additional information between the observation times. To
elucidate this point on a more intuitive level, we note that the
length of the time elapsed between $\tau_{k}$ and $\tau_{k+1}$\
carries additional information about the state of $\theta_{t}$
after $\tau_{k}.$ Specifically, if the frequency of observations
is proportional to the stock's volatility $v(\theta_{t})$,
$t\in[\hskip-.015in[\tau_{k},\tau_{k+1}]\hskip-.015in]$ ,
the larger values of $t-\tau_{k}$ might indicate lower values of
$v(\theta_{t})$.}
\end{remark}

\subsection{Volatility process}

A more precise description of the volatility process is in order
now. Let $(\mathbb{R},\mathcal{B}(\mathbb{R}))$ and
$(\mathbb{R}_{+}\times
\mathbb{R},\mathcal{B}(\mathbb{R}_{+})\otimes\mathcal{B}(\mathbb{R}))$
be measurable spaces with Borel $\sigma$-algebras. The volatility
process $\theta=(\theta_{t})_{t\geq0}$ is defined by the It\^{o}
equation
\begin{equation}
d\theta_{t}=b(t,\theta_{t})dt+\sigma(t,\theta_{t})dW_{t}+\int_{\mathbb{R}
}u(\theta_{t-},x)(\mu^{\theta}-\nu^{\theta})(dt,dx),
\label{2.1}
\end{equation}
where $W_{t}$ is a standard Wiener process and
$\mu^{\theta}=\mu^{\theta}(dt,dx)$ is a Poisson measure on $\big(
\mathbb{R}_{+}\times\mathbb{R},\mathcal{B}(\mathbb{R}_{+})
\otimes\mathcal{B}(\mathbb{R})\big)  $ with the compensator
$\nu^{\theta}(dt,dx)=K(dx)dt$, where $K(dx)$ is a $\sigma-$finite
non-negative measure on $\big(
\mathbb{R},\mathcal{B}(\mathbb{R})\big)  $. We assume that
$E\theta_{0}^{2}<\infty$, the functions $b(t,z),\sigma(t,z),$ and
$u(z,x)$ are Lipschitz continuous in $z$ uniformly with respect to
other variables, and
\[
|b(t,z)|+|\sigma(t,z)|^{2}+\int_{\mathbb{R}}|u(z,x)|^{2}K(dx)\leq
C(1+|z|^{2}).
\]
It is well known that under these assumptions \eqref{2.1}
possesses a unique strong solution adapted to $\mathbf{F}$, and
$E\theta_{t}^{2}<\infty$ for any $t\geq0$.

The generator $\mathcal{L}$ of the volatility process is given by
\begin{multline}
\mathcal{L}f(x):=b(t,x)f^{\prime}(x)+\frac{1}{2}\sigma^{2}(t,x)f^{\prime
\prime}(x)\label{2.10a}\\
+\int_{\mathbb{R}}\Big(f(x+u(x,y))-f(x)-f^{\prime}
(x)u(x,y)\Big)K(dy).\nonumber
\end{multline}

Before proceeding with the assumptions and main results we shall
introduce additional notation. Set
\begin{equation*}
m(s,t)=\int_{s}^{t}\left(
r(\theta_{u})-\frac{1}{2}v^{2}(\theta_{u})\right)du,
\end{equation*}
and
\begin{equation*}
\sigma^{2}(s,t)=\int_{s}^{t}v^{2}(\theta_{u})du.
\end{equation*}
For simplicity, it is assumed that $v^{2}(s,t)$ is bounded away
from zero. Let us denote by $\rho_{s,t}(y)$ the density function
of the normal distribution with mean $m(s,t)$ and the variance
$\sigma^{2}(s,t)$:
\begin{equation}
\rho_{s,t}(y):=\frac{1}{\sqrt{2\pi}\sigma(s,t)}e^{-\frac{(y-m(s,t))^{2}
}{2\sigma^{2}(s,t)}}.
\label{rho}
\end{equation}
Clearly, $\rho$ is the conditional density of the stock's
log-increments $X_{t}-X_{s}$ given $\theta$.

Let $\mathcal{F}^{\theta}=(\mathcal{F}_{t}^{\theta})_{t\geq0}$ be
the right-continuous filtration generated by
$(\theta_{t})_{t\geq0}$ and augmented by $\mathsf{P}$-zero sets
from $\mathcal{F}$. Denote by $G_{k}^{\theta}$ the conditional
distribution of $\tau_{k+1}$with respect to\footnote{Here and
below $\mathcal{F}^{1}\vee\mathcal{F}^{2}$ stands for the
$\sigma$-algebra generated by the $\sigma-$algebras
$\mathcal{F}^{1}\ $and $\mathcal{F}^{2}.$}
$\mathcal{F}^{\theta}\vee\mathcal{G}(k).$
That is,
$G_{k}^{\theta}$ is the distribution of the time of the next
observation, given previous history, and given $\theta$:
\begin{equation}
G_{k}^{\theta}(dt)  =\mathsf{P}\big(\tau_{k+1}\in
dt|\mathcal{F}^{\theta}\vee\mathcal{G}(k)\big).
\label{Gktheta}
\end{equation}
Without loss of generality we can and will assume that
$G_{k}^{\theta}(dt)  $ is the regular version of the
right hand side of (\ref{Gktheta}).

Let $N=(N_{t})_{t\geq0}$ be the counting process with interarrival
times: $\tau_0=0$, $(\tau_{k}-\tau_{k-1})_{k\geq1},$ that is
\begin{equation*}
N_{t}=\sum_{k\geq1}I(\tau_{k}\leq t).
\end{equation*}

\subsection{Assumptions}

The following assumptions will be in force throughout the paper:

A.0: For every $\mathcal{G}$-predictable and a.s. finite
stopping time $S$,
$$
\mathsf{P}(N_{S}-N_{S-}\neq0 |\mathcal{G}_{S-})=0 \ {\rm or} \ 1.
$$

\medskip
A.1: The Brownian motion $B$ is independent of
$\big(\theta,N\big)$.

\medskip A.2:
For every $k$, there exists a
$\mathcal{G}(k)$-measurable integrable random measure
$\Phi_{k}$ on $\mathcal{B}(\mathbb{R}_{+})$, so that,
for almost all $\omega \in\Omega,$ $\Phi_{k}\big([
0,\tau_{k}(\omega)]\big) =0$ and
$G_{k}^{\theta}$ is absolutely continuous with respect to
$\Phi_{k}$ $.$

Denote by $\phi(\tau_{k},t)=\phi(\theta,\tau_{k},t)$ the Radon-Nikodym derivative of
$G_{k}^{\theta}(dt)$ with respect to $\Phi_{k}(dt),$ i.e. for almost every $\omega,$
\begin{equation}
\phi(\tau_{k},t):=\frac{dG_{k}^{\theta}\big((\tau _{k},t]\big)}
{d\Phi_{k}\big((\tau_{k},t]\big)}.
\label{fika}
\end{equation}

Assumption A.0 is not essential for the derivation of the
filter. However, under this assumption the structure of the
optimal filter is simpler, and in the practical examples important
for this paper, this assumption holds anyway. In particular, A.0\
is verified if the conditional distribution
$
G_{k}^{\theta}=\mathsf{P}\big(\tau_{k+1}\leq
t|\mathcal{F}^{\theta} \vee\mathcal{G}(k)\big)
$
is absolutely continuous with respect to the Lebesgue
measure\footnote{More generally, it holds if the compensator of
the counting process $N_{t}$ is a continuous process.} or if the
arrival times $\tau_{k}$ are non-random.

The following two simple but important examples illustrate the
assumption A.2.

\begin{example}
\label{ex:cox copy(1)} Let ($\tau_{k})_{k\geq0}$ be \ the
jump times of a doubly stochastic Poisson process \ (Cox process)
with the intensity $n(\theta_{t}).$ In this case,
\[
\mathsf{P}\big(\tau_{k+1}\leq
t|\mathcal{F}^{\theta}\vee\mathcal{G}(k)\big )=
\begin{cases}
1-e^{-\int_{\tau_{k}}^{t}n(\theta_{s})ds} & ,t\geq\tau_{k}\\
0 & ,\text{otherwise}.
\end{cases}
\]
Then, one can take $\Phi_{k}(ds)=ds$ and
$\phi(\tau _{k},s)=n(\theta_{t})\exp\big(-\int_{\tau_{k}}^{s}n(\theta _{u})du\big)$.
If
$n(\theta_{t})=n$ is a constant, one could also choose
$$
\Phi_{k}(ds)=n\exp\big\{n(\tau_{k}-s)\big\}ds \quad\text{and}\quad
\phi(\tau_{k},s)=1.
$$
\end{example}

\begin{example}
\label{ex:constep0}
If the filtering is based on non-random
observation times $\tau_{k}$ (e.g., $\tau_{k}=kh$  where
$h$ is a fixed time step) then a natural choice would be
$$
\Phi_{k}(ds)=\delta_{\{\tau_{k+1}\}}(s)ds \ \text{and} \ \phi (\tau_{k},s)=1.
$$
\end{example}

For practical purposes, $\Phi_{k}(ds)$ must be
known or easily computable as soon as the the observations $(
\tau_{i},X_{i})_{i\leq k}$ become available. In contrast,
the Radon-Nikodym density $\phi(\tau_{k})  $ is, in
general, a function of the volatility process and is subject to
estimation.

We note that A.2 could be weakened slightly by replacing
$G_{k}^{\theta}$ by a regular version of the conditional \
distribution of $\tau_{k+1}$with respect to
$\mathcal{F}_{\tau_{k+1}-}^{\theta}\vee\mathcal{G}(k).$
The latter assumption would make the proof a little bit more
involved and we leave it to the interested reader.

\section{\bf Main results and introductory examples}

\label{Sec-3}

\subsection{Main result}

\label{sec-3.1ab} For a measurable function $f$ on $\mathbb{R}$
with $E|f(\theta_{t})|<\infty,$ define the conditional expectation
estimator $\pi_{t}(f)$ by
\begin{equation}
\pi_{t}(f):=E\big(f(\theta_{t})|\mathcal{G}_{t}\big)=\int_{\mathbb{R}}
f(z)\pi_{t}(dz),
\label{eq:p0}
\end{equation}
where $\pi_{t}(dz):=d\mathsf{P}(\theta_{t}\leq z|\mathcal{G}_{t})$
is the filtering distribution. (Note that we omit the argument
$\theta_{t}$ of $f$ in the estimator $\pi_{t}(f)$). In the spirit
of the Bayesian approach, it is assumed that the a priori
distribution
$$
\pi_{0}(dx)=\mathsf{P}\big(\theta_{0}\in dx\big)
$$
is given.

\medskip Let $\sigma\{\theta_{\tau_{k}}\}$ be the $\sigma$-algebra generated
by $\theta_{\tau_{k}}$. For $t>\tau_{k}$, let us define the
following \textit{structure functions}:

\begin{equation}
\psi_{k}(f;t,y,\theta_{\tau_{k}}):=E\Big(f(\theta_{t})\rho_{{\tau_{k},t}
}(y-X_{k})\phi(\tau_{k},t)\big|\sigma\big\{\theta_{\tau_{k}}\big\}\vee
\mathcal{G}(k)\Big),
\label{nado}
\end{equation}
and its integral with respect to $y$
\begin{equation}
\overline{\psi}_{k}(f;t,\theta_{\tau_{k}}):=\int_{\mathbb{R}}\psi_{k}\big(
f;t,y,\theta_{\tau_{k}}\big)
dy=E\Big(f(\theta_{t})\phi(\tau_{k}
,t)\big|\sigma\big\{\theta_{\tau_{k}}\big\}\vee\mathcal{G}(k)\Big),
\label{nado1}
\end{equation}
where $\rho$ and $\phi$ are given by (\ref{rho}) and (\ref{fika}),
respectively.\

If $f\equiv1$, the argument $f$\ in $\psi$ and \ $\bar{\psi}$ is
replaced by $1.$

Write
\[
\Phi_{k}(\{\tau_{k+1}\}):=\int_{0}^{\infty}I(t=\tau_{k+1})\Phi_{k}(dt),
\]
i.e. $\Phi_{k}(\{\tau_{k+1}\})$ is the jump of $\Phi_{k}(dt)$ at
$\tau_{k+1}$.

Finally, for $t\geq\tau_{k}$\ and a bounded function $f$, define
\[
\mathcal{M}_{k}(f;t,\pi_{t})
:=\frac{\pi_{\tau_{k}}( \bar{\psi}_{k}(f;t))  -\pi_{t-}(f)\pi_{\tau_{k}}
(\bar{\psi}_{k}(1;t))}{\int_{t}^{\infty}\pi_{\tau_{k} }(
\bar{\psi}_{k}(1;s))\Phi_{k}(ds)},
\]
whenever the denominator is not zero, and $\mathcal{M}_{k}(f;t,\pi_{t})=0$ if the
denominator is zero.

The main result of this paper is as follows:

\begin{theorem}
\label{mainthm} Assume A.0-A.2. Then for every measurable bounded
function $f$ in the domain of the generator $\mathcal{L}$ such
that $\int_{0} ^{t}E|\mathcal{L}f(\theta_{s})|ds<\infty$ for any
$t\geq0,$ the following system of equations holds{\rm :}

{\rm 1)} For every $k=0,1\ldots,$

\begin{equation}
\pi_{\tau_{k+1}}(f)=\frac{\pi_{\tau_{k}}(\psi_{k}(f;t,y))}{\pi_{\tau_{k}}
(\psi_{k}(1;t,y))}_{\big\{\substack{t=\tau_{k+1}\\y=X_{k+1}}
\big\}}-\mathcal{M}_{k}(f;t,\pi_{t})
_{_{\{t=\tau_{k+1}\}}}\cdot\Phi(\{\tau_{k+1}\}).
\label{eq:jump}
\end{equation}

{\rm 2)} For every $k=0,1\ldots$ and
$t\in]\hskip-1.5pt]\tau_{k},\tau _{k+1}[\hskip-1.5pt[$,
\begin{equation}
d\pi_{t}(f)=\pi_{t}(\mathcal{L}f)dt-\mathcal{M}_{k}(f;t,\pi_{t})\Phi_{k}(dt).
\label{eq:cont}
\end{equation}
\end{theorem}

\vspace{3mm}

\subsection{\textbf{Remarks}}
\mbox{} \vspace{1mm}

- 1. Equations (\ref{eq:jump}), (\ref{eq:cont}) form a closed
system of equations for the filter $\pi_{t}(f)$. It is often
convenient and customary (see e.g. \cite{Roz1}, \cite{Roz2} and
the references therein) to write a differential equation for a
measure-valued process $H_{t}(dx)$ in its
variational form, i.e. as the related system of equations for
$H_{t}(f)$ for all $f$ from a
sufficiently rich class of test functions belonging to the domain
of the operator $\mathcal{L}.$ \ In our setting, such a reduction
to the variational form is a necessity, since in some cases the
filtering measure $\pi_{s}(dx)=\mathsf{P}(\theta_{s}\in dx|\mathcal{G}_{s})
$may not belong to the
domain of $\mathcal{L}$. However, in the important examples
discussed below, there is no need to resort to the variational
form. The interested reader who is unaccustomed to the variational
approach might benefit from looking first into the examples at the
end of this section and in Section \ref{Sec-Ex}, where the
filtering equations are written as equations for posterior
distributions.

- 2. The system (\ref{eq:jump}) simplifies considerably if
\begin{equation}\label{m0}
\mathcal{M}_{k}(f;t,\pi_{t})
_{_{\{t=\tau_{k+1}\}}}\cdot\Phi(\{\tau_{k+1}\})=0, \ \text{\rm for all }k.
\end{equation}
Obviously, (\ref{m0}) holds if for all $k,$ $\Phi_{k}(dt)$
is continuous at $t=\tau_{k+1}$ as in the case when
$N_{t}$ is a Cox process. In fact, (\ref{m0}) holds true in many
other interesting cases, even when $\Phi_{k}(dt)$
has jumps at all $\tau_{k+1}$, as in the case of fixed observation
intervals (see Example \ref{ex:constep} below).\ We note then that
the following \textit{separation principle} holds.

\begin{corollary}
\label{cor:sp} Assume {\rm (\ref{m0})}. Then the filtering at the
observation times $\{\tau_{k}\}_{k\geq1}$ does not require
filtering between them; it is done by the Bayes type recursion{\rm
:}
\begin{equation*}
\pi_{\tau_{k+1}}(f)=\frac{\pi_{\tau_{k}}(\psi_{k}(f;t,y))}{\pi_{\tau_{k}}
(\psi_{k}(1;t,y))}_{\big\{\substack{t=\tau_{k+1}\\y=X_{k+1}}\big\}}.
\end{equation*}

\end{corollary}

- 3. Note that for high-frequency observations, even if condition
(\ref{m0})\ is not met, for all practical purposes, it may suffice
to compute the volatility estimates only at the observation times.
In that case, one would only use the relatively simple recursion
formula (\ref{eq:jump}), and disregard equation (\ref{eq:cont}).

- 4. Clearly, the \textquotedblleft structure functions" $\psi$
and $\bar {\psi}$ are of paramount importance for computing the
posterior distribution of the volatility process. We would like to
stress that these do not involve the observations and could be
pre-computed \textquotedblleft off-line" using just the \textit{a
priori} distribution. Then, \textquotedblleft on-line", when the
observations become available, one needs only to plug in the
obtained measurements $(\tau_{k},X_{k}),$ and to compute
$\pi_{t}(f)$ by recursion$.$ This feature is important for
developing efficient numerical algorithms.

- 5. Note also that for almost every $\omega\in\Omega,$ filtering
equation (\ref{eq:cont}) is a \textit{linear deterministic}
equation of Kolmogorov's type, rather than a \textit{nonlinear
stochastic} partial differential equation. The latter is typical
of the nonlinear filtering of diffusion processes. The
well-posedness and the regularity properties of equation
(\ref{eq:cont}) are well researched in the literature on second
order parabolic deterministic integro-differential equations (see
e.g. \cite{LM}, \cite{MP}, \cite{SK} and the references therein).

\begin{example}\label{examplechain}
(\textit{Volatility as a Markov Chain}.)
Let us now assume that the counting process is a Cox
process with intensity $n(\theta_{t})$, and take
$\phi(\tau_{k},s)=n(\theta_{t})e^{-\int_{\tau_{k}
}^{s}n(\theta_{u})du}$ and $\Phi_{k}(ds)=ds.$ Also
assume $\theta=(\theta_{t})_{t\leq T}$ is a homogeneous Markov
jump process taking values in the finite alphabet
$\{a_{1},\ldots,a_{M}\}$ with the intensity matrix
$\Lambda=\|\lambda(a_{i},a_{j})\|$ and the initial
distribution $p_{q}=\mathsf{P}(\theta_{0}=a_{q}),\ q=1,\ldots,M$.
(This is one of the two models of the state process discussed in
\cite{FR}.) In this case,
\[
\mathcal{L}f(\theta_{s})=\sum_{j}\lambda(\theta _{s},a_{j})f(a_{j}).
\]
Denote by $\theta_{t\text{ }}^{j}$ the process $\theta_{t\text{
}}$ starting from $a_{j}$, and
\begin{gather*}
p_{ji}(t):=\mathsf{P}\big(\theta_{t}=a_{i}|\theta_{0}=a_{j}\big), \quad
\pi_{j}(t)=\mathsf{P}\big(
\theta_{t}=a_{j}\big|\mathcal{G}_{t}\big),
\\
r_{ji}\left(  t,z\right)
:=E\big(e^{-\int_{0}^{t}n(\theta_{u}^{j})du}
\rho_{_{0,t}}^{j}(z)|\theta_{t}^{j}=a_{i}\big),
\end{gather*}
where $\rho_{_{0,t}}^{j}(z)$ is obtained by substituting
$\theta_{s}^{j}$ for $\theta_{s}$ in
$\rho_{_{0,t}}(z).$ It follows from Theorem \ref{mainthm} (for
details see Example \ref{exs1} ), with $f(\theta_{t})
:=I_{\{\theta_{t}=a_{i}\}},$ that
\begin{equation}
\pi_{i}(\tau_{k})=\frac{n(a_{i})
\sum_{j}r_{ji}(\tau _{k}-\tau_{k-1},X_{k}-X_{k-1})
p_{ji}(\tau_{k}-\tau _{k-1})
\pi_{j}(\tau_{k-1})}{\sum_{i,j}n(a_{i})r_{ji}(
\tau_{k}-\tau_{k-1},X_{k}-X_{k-1})p_{ji}(\tau
_{k}-\tau_{k-1})\pi_{j}(\tau_{k-1})}.
\label{chainTk}
\end{equation}
This recursion can be easily computed, once one computes
off-line the values $r_{ij}$. This example is
also treated in more detail in Section \ref{Sec-Ex}.
\end{example}

\section{\bf Proofs}
\label{Section3}

In the proof of the main result we want to show that
\[
d\pi_{t}(f) =\pi_{t}(\mathcal{L}f)dt +dM_{t},
\]
where $M_{t}$ is a martingale, and then we find a (integral)
martingale representation of $M_{t}$ with respect to the measure
$\mu-\nu$, where $\nu$ is a compensator of $\mu$. We first find
the compensator.

\subsection{$\b{(\mathcal{G}_{t})}$-compensator of $\b{\mu}$}

\label{subsec-3} Denote by $\mathcal{P}(\mathcal{G})$ be the
predictable $\sigma$-algebra on $\Omega\times[0,\infty)$
with respect to $\mathcal{G}$ and and set
\[
\widetilde{\mathcal{P}}(\mathcal{G})=\mathcal{P}(\mathcal{G})\otimes
\mathcal{B}(\mathbb{R}).
\]

A nonnegative random measure $\nu(dt,dy)$ on
$\widetilde{\mathcal{P} }(\mathcal{G})$ is called a
$\widetilde{\mathcal{P}}(\mathcal{G})$ -compensator of $\mu$ if
for any $\widetilde{\mathcal{P}}(\mathcal{G})$-measurable,
nonnegative function $\varphi(t,y)=\varphi(\omega,t,y)$,
\begin{equation*}
\begin{split}
\mathrm{(i)}\quad &
\int_{0}^{t}\int_{\mathbb{R}}\varphi(s,y)\nu(ds,dy)
\ \text{is ${\mathcal{P}}(\mathcal{G})$-measurable}\\
\mathrm{(ii)}\quad &  E\int_{0}^{\infty}\int_{\mathbb{R}}\varphi
(t,y)\mu(dt,dy)=E\int_{0}^{\infty}\int_{
\mathbb{R}}\varphi(t,y)\nu(dt,dy).
\end{split}
\end{equation*}
Let
$
G_{k}(ds,dx)=G_{k}(\omega,ds,dx)
$
be a regular version of the conditional distribution of
$
\big(
\tau_{k+1} ,X_{k+1}\big)  $ given $\mathcal{G}\left(  k\right)
$
( it is assumed that $G_{k}([0,\tau_{k}],dx)=0$):
\begin{align}\label{GGG}
\mathsf{G}_{k}(dt,dy)=d\mathsf{P}\big(\tau_{k+1}\le
t,X_{k+1}\le y|\mathcal{G}(k)\big).
\end{align}
Denote $G_{k}(ds)=G_{k}(dt,\mathbb{R}),$ that is, $G_{k}(t)=\mathsf{P}(\tau_{k+1}\leq
t~|~\mathcal{G}(k))$ (with probability one).

By Theorem III.1.33 \cite{JS} \ (see also Proposition 3.4.1 in
\cite{LSMar}),
\begin{equation*}
{\nu}(dt,dy)=\sum_{k\geq0}I_{]\hskip-1.5pt]\tau_{k},\tau_{k+1}]\hskip-1.5pt]}
(t)\frac{G_{k}(dt,dy)}{G_{k}([t,\infty),\mathbb{R})},
\end{equation*}
We now derive a representation, suitable for the filtering
purposes, of the
$\widetilde{\mathcal{P}}(\mathcal{G})$-compensator $\nu$ in terms
of the structure functions (\ref{nado}), (\ref{nado1}), and the
posterior distribution of $\theta\ $.

\begin{lemma}
\label{lem-4.1} The
$\widetilde{\mathcal{P}}(\mathcal{G})-$compensator $\nu$ admits
the following version{\rm :}
\begin{equation}
\nu(dt,dy)=\sum_{k\geq0}I_{]\hskip-1.5pt]\tau_{k},\tau_{k+1}]\hskip-1.5pt]}
(t)\frac{\pi_{\tau_{k}}(\psi_{k}(1;t,y))}{\int_{t}^{\infty}\pi_{\tau_{k}
}(\overline{\psi}_{k}(1;s))\Phi_{k}(ds)}\Phi_{k}(dt)dy.
\label{eq:nu}
\end{equation}

\end{lemma}

\begin{proof}
By A.1 for $t>\tau_{k}$, with probability 1,

\begin{equation*}
\begin{split}
&  \mathsf{P}\big(\tau_{k+1}\leq t,X_{k+1}\leq
y|\mathcal{F}^{\theta}\vee\mathcal{G}(k)\big)
\\
&=E\Big(\mathsf{P}\big(\tau_{k+1}\leq t,X_{k+1}\leq
y|\mathcal{F}^{\theta }\vee\mathcal{G}(k)
\vee\sigma(\tau_{k+1})
\big)\big|\mathcal{F}^{\theta}\vee\mathcal{G}(k)\Big)
\\
&  =E\Big(I_{(\tau_{k+1}\leq t)}\mathsf{P}\big(
X_{k+1}\leq y|\mathcal{F}^{\theta}\vee\mathcal{G}\left(  k\right)
\vee \sigma\left(  \tau_{k+1}\right)  \big)
|\mathcal{F}^{\theta}\vee\mathcal{G}(k)\Big)
\\
&=E\left(I_{(\tau_{k+1}\leq t)}\int_{-\infty}^{y}\rho _{\tau_{k},\tau_{k+1}}(
z-X_{k})dz|\mathcal{F}^{\theta}
\vee\mathcal{G}(k)\right)
\\
& =\int_{\tau_{k}}^{t}\int_{-\infty}^{y}\rho_{\tau_{k},s}(z-X_{k}
)dzG_{k}^{\theta}(ds)  ,
\end{split}
\end{equation*}
where we recall that $G_{k}^{\theta}$ is a regular version of the
conditional \ distribution of $\tau_{k+1}$with respect to
$\mathcal{F}^{\theta} \vee\mathcal{G}(k).$ Thus, by
A.2, for $t>\tau_{k}$, with probability 1,

\begin{gather}\label{eq:sep}
  \mathsf{P}\big(\tau_{k+1}\leq t,X_{k+1}\leq
y|\mathcal{F}^{\theta}
\vee\mathcal{G}(k)\big)
\nonumber\\
  =\int_{\tau_{k}}^{t}\int_{-\infty}^{y}\rho_{\tau_{k},s}(
z-X_{k})\phi(\tau_{k},s) dz\Phi_{k}(ds).
\end{gather}
By (\ref{nado}), using notation (\ref{eq:p0}), we
see that
\[
E\big(E\big[\phi(\tau_{k},s)\rho_{\tau_{k},s}(z-X_{k})|\sigma\{\theta
_{\tau_{k}}\}\vee\mathcal{G}(k)
\big]|\mathcal{G}(k)
\big)=\pi_{\tau_{k}}(\psi_{k}(1;s,z)).
\]
This, together with (\ref{eq:sep}), yields, recalling definition
(\ref{GGG}),
\begin{equation*}
G_{k}\big(ds,dz)=\pi_{\tau_{k}}(\psi_{k}(1;s,z))\Phi_{k}(ds)dz.
\end{equation*}
In the same way, for $t>\tau_{k}$, with probability 1,

\begin{equation}
G_{k}\big([t,\infty],\mathbb{R})=\int_{t}^{\infty}\pi_{\tau_{k}}
(\overline{\psi}_{k}(1;s))\Phi_{k}(ds).
\label{eq:denom}
\end{equation}
This completes the proof.
\end{proof}

\begin{remark}
\label{re:den0}
{\rm If the right hand of \ (\ref{eq:denom}) is zero,
then
$
\mathsf{P}\big(\tau_{k+1}\geq t|\mathcal{G}(k)\big)=0.
$
Hence,
$I_{]\hskip-1.5pt]\tau_{k},\tau_{k+1}]\hskip-1.5pt]}(t)=0$ with
probability 1 and, by the $\ 0/0=0$ convention, the corresponding
term in (\ref{eq:nu}) is zero.}
\end{remark}

\subsection{Semimartingale representation of the optimal filter}

\label{sec-4} In this section we will prove the following result.

\begin{theorem}
\label{theo-4.1} For any bounded function $f$ from the domain of
the operator $\mathcal{L}$ such that
$\int_{0}^{t}E|\mathcal{L}f(\theta_{s})|ds<\infty$ for all
$t<\infty$, the differential of the optimal filter $\pi_{s}(f)$ is
given by equation
\begin{align}\label{filteq}
d\pi_{s}(f)  &  =\pi_{s}(\mathcal{L}f)ds
\\
&
+\int_{\mathbb{R}}\Big(\sum_{k\geq0}I_{]\hskip-1.5pt]\tau_{k},\tau
_{k+1}]\hskip-1.5pt]}(s)\frac{\pi_{\tau_{k}}(\psi_{k}(f;s,y))}{\pi_{\tau_{k}
}(\psi_{k}(1;s,y))}-\pi_{s-}(f)\Big)(\mu-\nu)(ds,dy).\nonumber
\end{align}

\end{theorem}

\begin{proof}
It suffices to verify the statement for twice continuously
differentiable functions $f$ with $f,f^{\prime}f^{\prime\prime}$
bounded. By It\^{o}'s formula,

\[
\begin{split}
&
f(\theta_{t})=f(\theta_{0})+\int_{0}^{t}\mathcal{L}f(\theta_{s})ds+\int
_{0}^{t}f^{\prime}(\theta_{s})\sigma(\theta_{s})dW_{s}\\
&
+\int_{0}^{t}\int_{\mathbb{R}}f^{\prime}(\theta_{s-})u(\theta_{s-}
,x)(\mu^{\theta}-\nu^{\theta})(ds,dx).
\end{split}
\]
Denote
\[
L_{t}=\int_{0}^{t}f^{\prime}(\theta_{s})\sigma(\theta_{s})dW_{s}+\int_{0}
^{t}\int_{\mathbb{R}}f^{\prime}(\theta_{s-})u(\theta_{s-},x)(\mu^{\theta}
-\nu^{\theta})(ds,dx).
\]
Then, we have
\begin{align*}
\pi_{t}(f) &  =E\big(f(\theta_{0})|\mathcal{G}_{t}\big)\\
&  +E\Bigg(\int_{0}^{t}\mathcal{L}f(\theta_{s})ds\Big|\mathcal{G}
_{t}\Bigg)+E\big(L_{t}|\mathcal{G}_{t}\big).
\end{align*}
Set
\[
\begin{split}
M_{t} &  =\big\{E\big(f(\theta_{0})|\mathcal{G}_{t}\big)-\pi_{0}(f)\big\}\\
&
\quad+\Bigg\{E\Bigg(\int_{0}^{t}\mathcal{L}f(\theta_{s})ds\Big|\mathcal{G}
_{t}\Bigg)-\int_{0}^{t}\pi_{s}\big(\mathcal{L}f\big)ds\Bigg\}+E\big(L_{t}
|\mathcal{G}_{t}\big).
\end{split}
\]
Obviously, the process
$E\big(f(\theta_{0})|\mathcal{G}_{t}\big)-\pi_{0}(f)$ is a
$\mathcal{G}_{t}$-martingale. Process $L_{t}$ is a
$\mathcal{F}_{t} $-martingale. Since
$\mathcal{G}_{t}\subseteq\mathcal{F}_{t}$, for $t>t^{\prime},$
\[
E\big(E(L_{t}|\mathcal{G}_{t})|\mathcal{G}_{t^{\prime}}\big)=E\big(E(L_{t}
|\mathcal{F}_{t^{\prime}})|\mathcal{G}_{t^{\prime}}\big)=E(L_{t^{\prime}
}|\mathcal{G}_{t^{\prime}}).
\]
Consequently, $E(L_{t}|\mathcal{G}_{t})$ is a martingale too.

Finally,
$E\big(\int_{0}^{t}\mathcal{L}f(\theta_{s})ds|\mathcal{G}_{t}\big)-\int
_{0}^{t}\pi_{s}\big((\mathcal{L}f)\big)ds$ is also a
$\mathcal{G}_{t} $-martingale.
Indeed, for $t>s>t^{\prime},$ we
have $E\big(\pi_{s}
\big(\mathcal{L}f)\big|\mathcal{G}_{t^{\prime}}\big)=E\big(\mathcal{L}
f(\theta_{s})|\mathcal{G}_{t^{\prime}}\big)$ which yields
\begin{align*}
&  E\left[
E\Bigg(\int_{0}^{t}\mathcal{L}f(\theta_{s})ds\Big|\mathcal{G}
_{t}\Bigg)-\int_{0}^{t}\pi_{s}(\mathcal{L}f)ds\Bigg|\mathcal{G}_{t^{\prime}
}\right]
\\
& \quad=E\Bigg(\int_{0}^{t^{\prime}}\mathcal{L}f(\theta_{s}
)ds\Big|\mathcal{G}_{t^{\prime}}\Bigg)-\int_{0}^{t^{\prime}}\pi_{s}
(\mathcal{L}f)ds.
\end{align*}

Thus, $M_{t}$ is a $\mathcal{G}_{t}$-martingale. In particular,
this means that $\pi_{t}(f)$ is a $\mathcal{G}$ - semimartingale
with paths in the Skorokhod space
$\mathbb{D}_{[0,\infty)}(\mathbb{R})$, so that $\pi_{t}(f)$ is a
right continuous process with limits from the left. By the
Martingale Representation Theorem ( see e.g. Theorem 1 and Problem
1.c in Ch.4, \S 8. in \cite{LSMar}),
\[
M_{t}=\int_{0}^{t}\int_{\mathbb{R}}H(s,y)(\mu-\nu)(ds,dy).
\]
It is a standard fact that
$
\mathsf{P}(N_{S}-N_{S-}\neq0
|\mathcal{G}_{S-})={\nu}(\{S\},\mathbb{R}_{+}).
$
Hence, due to assumption A.0, by Theorem 4.10.1 from \cite{LSMar}
(see formulae (10.6) and (10.15)),
\begin{equation*}
H(t,y)=\mathsf{M}_{\mu}^{\mathsf{P}}\big(\triangle
M|\widetilde{\mathcal{P} }(\mathcal{G})\big)(t,y),
\end{equation*}
where $\triangle M_{t}=M_{t}-M_{t-}$ and the conditional
expectation
$\mathsf{M}_{\mu}^{\mathsf{P}}\big(g|\widetilde{\mathcal{P}}(\mathcal{G}
)\big)$ is defined by the following relation (see, e.g.
\cite{LSMar}, Ch. 2, \S 2 and Ch. 10, \S 1): for any
$\widetilde{\mathcal{P}}(\mathcal{G} )$-measurable bounded and
compactly supported function $\varphi(t,y),$
\begin{equation*}
E\int_{0}^{\infty}\int_{\mathbb{R}}\varphi(t,y)g_{t}\mu(dt,dy)
=E\int_{0}^{\infty}\int_{\mathbb{R}}\varphi(t,y)\mathsf{M}_{\mu}^{\mathsf{P}
}\big(g\big|\widetilde{\mathcal{P}}(\mathcal{G})\big)(t,y)\nu(dt,dy).
\end{equation*}
By Lemma 4.10.2, \cite{LSMar},
\begin{equation}
\mathsf{M}_{\mu}^{P}\big(\pi_{t}(f)\big|\widetilde{\mathcal{P}}(\mathcal{G}
)\big)(t,y)=\mathsf{M}_{\mu}^{P}\big(f\big|\widetilde{\mathcal{P}}
(\mathcal{G})\big)(t,y).
\label{4.10a}
\end{equation}
Since, $\pi_{t-}(f)$ is
$\widetilde{\mathcal{P}}(\mathcal{G})$-measurable (which implies
$\mathsf{M}_{\mu}^{\mathsf{P}}(\pi_{-}(f)|\widetilde
{\mathcal{P}}(\mathcal{G}))(t,y)=\pi_{t-}(f)$ ), by
(\ref{4.10a}),
\begin{gather*}
\mathsf{M}_{\mu}^{\mathsf{P}}\big(\triangle
M\big|\widetilde{\mathcal{P}
}(\mathcal{G})\big)(t,y)
=\mathsf{M}_{\mu}^{\mathsf{P}}\big(\pi_{t}(f)-\pi_{t-}(f)\big|\widetilde
{\mathcal{P}}(\mathcal{G})\big)(t,y)\nonumber\\
\\
=\mathsf{M}_{\mu}^{\mathsf{P}}\big(f\big|\widetilde{\mathcal{P}
}(\mathcal{G})\big)(t,y)-\pi_{t-}(f).
\end{gather*}
To complete the proof one needs to show that
\begin{equation}
\mathsf{M}_{\mu}^{\mathsf{P}}\big(f(\theta_{.})
\big|\widetilde
{\mathcal{P}}(\mathcal{G})\big)(s,y)=\sum_{k\geq0}I_{]\hskip-1.5pt]\tau
_{k},\tau_{k+1}]\hskip-1.5pt]}(s)\frac{\pi_{\tau_{k}}(\psi_{k}(f;s,y))}
{\pi_{\tau_{k}}(\psi_{k}(1;s,y))}.
\label{gg}
\end{equation}
To prove (\ref{gg}), it suffices to demonstrate that for any
$\widetilde {\mathcal{P}}(\mathcal{G})$-measurable bounded and
compactly supported function $\varphi(t,y),$
\begin{gather}
E\sum_{k\geq0}\int_{(\tau_{k},\tau_{k+1}]\cap(\tau_{k},\infty)}\int
_{\mathbb{R}}\varphi(t,y)\frac{\pi_{\tau_{k}}(\psi_{k}(f;t,y))}{\pi_{\tau_{k}
}(\psi_{k}(1;t,y))}\nu(dt,dy)
\nonumber\\
=E\int_{0}^{\infty}\int_{\mathbb{R}}
\varphi(t,y)f(\theta_{t})\mu(dt,dy).
\label{show}
\end{gather}

By monotone class arguments, we can assume that $\varphi(t,x)
=v(t)g(x)$, where
$v(t)$ is a $\mathcal{P}(\mathcal{G})$-measurable
process and $g(x)$ is a continuous function on
$\mathbb{R}$. By Lemma III.1.39 \cite{JS}, since ${v} ${$(t)$ is
$\mathcal{P}(\mathcal{G})-$measurable, it must be
of the form
\begin{equation}
v(t)=v_{0}+\sum_{k\geq1}^{\infty}v_{k}(t)I_{]\hskip-1.5pt]\tau_{k},\tau_{k+1}]
\hskip-1.5pt]}(t),
\label{4.11b}
\end{equation}
where $v_{0}$ is a constant and $v_{k}\left(  t\right)  $ are
$\mathcal{G} \left(  k\right)  \otimes\mathcal{B}(\mathbb{R}_{+})$-measurable functions. }

Owing to (\ref{4.11b}) and Lemma \ref{lem-4.1}, in order to prove
(\ref{show}), it suffices to verify the equality
\begin{gather}
E\left[
\int_{(\tau_{k},\tau_{k+1}]\cap(\tau_{k},\infty)}\int_{\mathbb{R}
}g(y)v_{k}(t)\frac{\pi_{\tau_{k}}(\psi_{k}(f;t,y))}{\pi_{\tau_{k}}(\psi
_{k}(1;t,y))}\Phi_{k}(dt)dy\right]
\nonumber\\
\nonumber\\
=E\left[
v_{k}(\tau_{k+1})g(X_{k+1})f(\theta_{\tau_{k+1}})1_{\{
\tau_{k+1}<\infty\}}\right],
\label{enough}
\end{gather}

The next step follows the ideas of Theorem III.1.33 \cite{JS}. We
have
\begin{align*}
&  E\left[
v_{k}(\tau_{k+1})g(X_{k+1})f(\theta_{\tau_{k+1}})1_{\{
\tau_{k+1}<\infty\}}\right]
\\
&  =E\left[E\big(
v_{k}(\tau_{k+1})g(X_{k+1})f(\theta_{\tau_{k+1} })1_{\{
\tau_{k+1}<\infty\}}|\mathcal{G}(k)
\vee\mathcal{F}^{\theta}\big)\right]
\\
&=E\left(
\int_{(\tau_{k},\infty)}\int_{\mathbb{R}}v_{k}(s)g(y)E\left[
f(\theta_{s})G_{k}^{\theta}\left(  ds,dy\right)
|\mathcal{G}\left(  k\right) \right]\right),
\end{align*}
where, as before, $G_{k}^{\theta}$ $(ds,dy)$ is a regular version
of the conditional \ distribution of $\big(
\tau_{k+1},X_{k+1}\big)$ with respect to
$\mathcal{F}^{\theta}\vee\mathcal{G}(k).$

By Fubini Theorem, and recalling notation (\ref{GGG}),
\begin{align*}
&  E\left(
\int_{(\tau_{k},\infty)}\int_{\mathbb{R}}v_{k}(s)g(y)E\left[
f(\theta_{s})G_{k}^{\theta}(ds,dy)
|\mathcal{G}(k)
\right]\right)
\\
&=E\left(
\int_{(\tau_{k},\infty)}\int_{\mathbb{R}}v_{k}(s)g(y)\frac
{E\left[f(\theta_{s})G_{k}^{\theta}(ds,dy)
|\mathcal{G}(k)\right]}{G_{k}\left(\left[s,\infty\right];\mathbb{R}\right)
}\int_{[s,\infty]}G_{k}(du,\mathbb{R})\right)
\\
&  =E\left(
\int_{\tau_{k}}^{\tau_{k+1}}\int_{\mathbb{R}}v_{k}(s)g(y)\frac
{E\left[f(\theta_{s})G_{k}^{\theta}(ds,dy)
|\mathcal{G}(k)\right]}{G_{k}\left(\left[
s,\infty\right];\mathbb{R}\right) }\right).
\end{align*}
By (\ref{eq:sep}),
\begin{equation*}
G_{k}^{\theta}(ds,dy)
=\rho_{{\tau_{k},s}}(z-X_{k})\phi(\tau _{k},s)\Phi_{k}(ds)dy.
\end{equation*}

Hence, for $s>\tau_{k}$,
\begin{align*}
&  E\left[f(\theta_{s})G_{k}^{\theta}(ds,dy)
|\mathcal{G}(k\right]
\\
&  =E\Big(E\left(f(\theta_{s})\rho_{{\tau_{k},s}}(y-X_{k})\phi(\tau
_{k},s)|\sigma\{\theta_{\tau_{k}}\}
\vee\mathcal{G}(k)\right)\big|\mathcal{G}(k)\Big)\Phi_{k}(ds)dy
\\
&=\pi_{\tau_{k}}(\psi_{k}(f;s,y))dy\Phi_{k}(ds).
\end{align*}

This, together with (\ref{eq:denom}), yields

\begin{align*}
&
E\left(\int_{\tau_{k}}^{\tau_{k+1}}\int_{\mathbb{R}}v_{k}(s)g(y)\frac
{E\left[f(\theta_{s})G_{k}^{\theta}(ds,dy)
|\mathcal{G}(k)\right]}{G_{k}\left(\left[
s,\infty\right];\mathbb{R}\right)}\right)
\\
&=E\left(
\int_{\tau_{k}}^{\tau_{k+1}}\int_{\mathbb{R}}v_{k}(s)g(y)\frac
{\pi_{\tau_{k}}(\psi_{k}(f;s,y))dy}{\int_{s}^{\infty}\pi_{\tau_{k}}
\left(\bar{\psi}(1;t)\right)\Phi_{k}(dt)}\Phi_{k}(ds)\right)  ,
\end{align*}
so that (\ref{enough}) is satisfied, and the proof
follows.
\end{proof}

\subsection{Proof of Theorem \ref{mainthm}}
\label{sec-6}
In this section we show that Theorem \ref{mainthm}
follows from Lemma \ref{lem-4.1} and Theorem \ref{theo-4.1}.

\begin{proof}
Firstly, we note that the stochastic integral in the right hand side of
(\ref{filteq}) can be written as the difference of the integrals
with respect to $\mu$ and $\nu.$ Indeed, since $f$ is bounded,
this follows from \cite{JS}, Proposition II.1.28.

By applying Lemma \ref{lem-4.1} and integrating over $y$ one gets
that for $t\in]\hskip-1.5pt]\tau_{k},\tau_{k+1}]\hskip-1.5pt],$

\begin{align*}
&
\int_{\mathbb{R\times}(\tau_{k},t]}\Big(\frac{\pi_{\tau_{k}}(\psi
_{k}(f;s,y))}{\pi_{\tau_{k}}(\psi_{k}(1;s,y))}-\pi_{s-}(f)\Big)\nu(ds,dy)\\
&  =\int_{(\tau_{k},t]}\frac{\pi_{\tau_{k}}\left(
\bar{\psi}_{k}(f;s)\right) -\pi_{s-}(f)\pi_{\tau_{k}}\left(
\bar{\psi}_{k}(1;s)\right)}{\int
_{s}^{\infty}\pi_{\tau_{_{k}}}\left(\bar{\psi}_{k}(1;u)\right)
\Phi _{k}(du)}\Phi_{k}(ds).
\end{align*}

This equation verifies that \eqref{eq:cont} follows from the
semimartingale representation \eqref{filteq}, for $t$ between the
consecutive observation times.

For the jump part \eqref{eq:jump}, we note that
\[
\int_{0}^{t}\int_{\mathbb{R}}\pi_{s-}(f)\mu(ds,dy)=
\sum_{\tau_{k+1}\leq t}\pi_{(
\tau_{k+1})-}(f)
\]
and
\[
\int_{0}^{t}\int_{\mathbb{R}}\frac{\pi_{\tau_{k}}(\psi_{k}(f;s,y))}{\pi
_{\tau_{k}}(\psi_{k}(1;s,y))}\mu(ds,dy)
=\sum_{\tau_{k+1}\leq
t}\frac{\pi_{\tau_{k}}(\psi_{k}(f;s,y))}{\pi_{\tau_{k}}(\psi_{k}(1;s,y))}
_{\big\{\substack{s=\tau_{k+1}\\y=X_{k+1}}\big\}}.
\]
Now, (\ref{filteq}) can be rewritten as follows:
\begin{equation}\label{closeq1}
\begin{aligned}
\pi_{t}(f)&=\pi_{0}(f)
+\int_{0}^{t}\pi
_{s}\big(\mathcal{L}f\big)ds
\\
&  +\sum_{\tau_{k+1}\leq t}\left(
\frac{\pi_{\tau_{k}}(\psi_{k}(f;s,y))}
{\pi_{\tau_{k}}(\psi_{k}(1;s,y))}
_{\big\{\substack{s=\tau_{k+1}
\\y=X_{k+1}}\big\}} -\pi_{(\tau_{k+1})-}(f)
\right)
\\
&
-\sum_{k\geq0}\int_{(\tau_{k},t\wedge\tau_{k+1}]}\mathcal{M}_{k}(f;s,\pi_{s})
\Phi_{k}(ds).
\end{aligned}
\end{equation}
Suppose $t\in]\hskip-1.5pt]\tau_{k},\tau_{k+1}[\hskip-1.5pt[.$
Then,
\begin{align*}
\pi_{t}(f)&=\pi_{\tau_{k}}(f)
\\
&+\int_{\tau_{k}}^{t}\pi_{s}\big(\mathcal{L}f)ds
-\int_{\tau_{k} }^{t}\mathcal{M}_{k}(f;s,\pi_{s})
\Phi_{k}(ds).
\end{align*}
It follows that
\begin{align*}
&  \pi_{(\tau_{k+1})-}(f)
\\
&=\pi_{\tau_{k}}(f)
+\int_{\tau_{k}}^{\tau_{k+1}}\pi _{s}\big(\mathcal{L}f\big)
ds-\int_{\tau_{k}}^{(\tau _{k+1})-}\mathcal{M}_{k}(f;s,\pi_{s})\Phi_{k}(ds).
\end{align*}
Therefore, from (\ref{closeq1}),
\[
\pi_{\tau_{k+1}}(f)
=\frac{\pi_{\tau_{k}}(\psi_{k}(f;s,y))}
{\pi_{\tau_{k}}(\psi_{k}(1;s,y))}
_{\big\{\substack{s=\tau_{k+1}
\\y=X_{k+1}}\big\}}
-\mathcal{M}_{k}(f;t,\pi_{t})
_{\{t=\tau_{k+1}\}} \Phi(\{\tau_{k+1}\}).
\]

This completes the proof.
\end{proof}

\section{\bf Examples}
\label{Sec-Ex}

In this Section we consider some important special cases of
Theorem \ref{mainthm}.

\begin{example}
\label{exs1}
(\textit{Markov chain volatility and Cox process
arrivals.}) Recall the setting of Example \ref{examplechain} and
its notation $r_{ij}$, $\pi_{j}(t),$ and $\theta^{j}$. It follows
from Example \ref{ex:cox copy(1)} that in this case
$\Phi_{k}(\{  \tau_{k+1}\})=0$ for all $k$'s.
Hence the second term in the right hand side of equation
(\ref{eq:jump}) is zero. By (\ref{nado}), for $f(\theta
_{t})=\mathbf{1}_{\{\theta_{t}=a_{i}\}}$ and $t>\tau_{k}$,
\[
\psi_{k}(f;t,y,\theta_{\tau_{k}})=n(a_{i})\left[
E\big(I_{\{\theta_{t}=a_{i}\}}e^{-\int_{s}^{t}n(\theta_{u}
)du}\rho_{_{s,t}}(y-x)|\theta_{s}\big)\right]_{\big\{\substack
{s=\tau_{k}\\x=X_{k}}\big\}}.
\]
Thus, owing to the homogeneity of $\theta_{t},$ for $t$
$>\tau_{k},$
$$
\begin{aligned}
&\qquad \pi_{\tau_{k}}(\psi_{k}(f;t,y))
\\
&=\sum_{j}n(a_{i})E\Big(I_{\{\theta_{t}=a_{i}\}}e^{-\int
_{s}^{t}n(\theta_{u})du}\rho_{s,t}(y-x)\big|\theta_{s}=a_{j}
\Big)_{\big\{\substack{s=\tau_{k}\\x=X_{k}}\big\}}\pi_{j}(\tau_{k})
\\
&=\sum_{j}n(a_{i})E\Big(I_{\{  \theta_{t-s}^{j}=a_{i}\}}
e^{-\int_{0}^{t-s}n(\theta_{u})du}\rho_{0,t-s}^{j}
(y-x)\Big)_{\big\{\substack{s=\tau_{k}\\x=X_{k}}\big\}}\pi_{j}(\tau_{k})
\\
&=\sum_{j}n(a_{i})E\Big[I_{\{\theta_{t-s}^{j}=a_{i}\}}
E\Big(e^{-\int_{0}^{t-s}n(\theta_{u})du}\rho_{0,t-s}^{j}(y-x)\big|\theta
_{t-s}^{j}\Big)\Big]_{\hskip
-.05in\big\{\substack{s=\tau_{k}\\x=X_{k}}\big\}}
\pi_{j}(\tau_{k})
\\
&=\sum_{j}n(a_{i})r_{ji}(t-\tau_{k},y-X_{k})
p_{ji}(t-\tau_{k})\pi_{j}(\tau_{k}).
\end{aligned}
$$
Similar formula holds for the denominator of the first term of the
right hand side of the equation. Now equation \eqref{chainTk} follows from
\eqref{eq:jump}.

Repeating the previous calculations and using the notation
\[
\bar{r}_{ji}(t)
:=E\big(e^{-\int_{0}^{t}n(\theta_{u}^{j})du}|\theta_{t}^{j}=a_{i}\big),
\]
it is readily checked that, for $t>\tau_{k}$,
\[
\pi_{\tau_{k}}\big(\bar{\psi}_{k}(\mathbf{1}_{\{\theta_{t}=a_{i}
\}};t)\big)=n(a_{i})
\sum_{j}\pi_{j}(\tau_{k})\bar{r}_{ji}(t-\tau_{k})p_{ji}(t-\tau_{k})
\]
and
\[
\pi_{\tau_{k}}\big(\bar{\psi}_{k}(1,t)\big)=\sum_{i,j}\pi_{j}(\tau
_{k})n(a_{i})\bar{r}_{ji}(t-\tau_{k})p_{ji}(t-\tau _{k})
\]
which are needed in computing \eqref{eq:cont}. It is easily
verified that in the setting of this example, equation
\eqref{eq:cont} reduces to the following:
\begin{equation}
d\pi_{i}(t)=\sum_{j}\lambda(a_{j},a_{i})
\pi_{j}(t)dt+\bar {D}(\tau_{k},t)
\pi_{i}(t)dt+D_{i}(\tau_{k},t)dt,
\label{eq:ko2}
\end{equation}
where
\begin{align*}
D_{i}(\tau_{k},t)&=-\frac{n(a_{i})
\sum _{j}\bar{r}_{ji}(t-\tau_{k})p_{ji}(t-\tau_{k})
\pi_{j}(\tau _{k})}{\int_{t}^{\infty}\sum_{i,j}n(a_{i})\bar{r}_{ji}
(s-\tau_{k})p_{ji}(s-\tau_{k})\pi_{j}(\tau_{k})ds}
\\
\bar{D}(\tau_{k},t)&=\frac{\sum_{l,j}n(a_{l})\bar{r}_{jl}(t-\tau_{k})p_{jl}
(t-\tau_{k})\pi_{j}(\tau_{k})}{\int_{t}^{\infty}\sum_{i,j}n(a_{i})
\bar{r}_{ji}(s-\tau _{k})p_{ji}(s-\tau_{k})
\pi_{j}(\tau_{k})ds}.
\end{align*}
Note that equation (\ref{eq:ko2})\ is considered for a fixed
$\omega$ and $t>\tau_{k}(\omega).$ Therefore,
$\tau_{k}$ and $\pi_{\cdot }(\tau_{k})$ should be viewed as known
quantities.
\end{example}

\begin{example}
\label{Poisson}(\textit{Poisson arrivals.})
Suppose that the
interarrival times between the observations are exponential  with
constant intensity $n(\theta)\equiv\lambda$. In other words,
$N_{t}$ is Poisson process with constant parameter $\lambda.$ In
this case, the volatility process $\theta$ is independent of
$N_{t}.$ Then, on the interval $\tau_{k}<t<\tau_{k+1},$
equation (\ref{eq:ko2}) reduces to
\begin{align}\label{eq:po}
d\pi_{i}(t)&=\sum_{j}\lambda(a_{j},a_{i})
\pi_{j}
(t)dt
\nonumber\\
&  -\lambda\Big(\sum_{j}p_{ji}(t-\tau_{k})
\pi_{j}(\tau_{k})-\pi _{i}(t)\Big)dt.
\end{align}
On the other hand, owing to the independence of $N$ and
$\theta,$ it is readily checked that on the interval
$\tau_{k}<t<\tau_{k+1},$
\[
\pi_{i}(t)={\ }\sum_{j}p_{ji}(t-\tau_{k})
\pi_{j}(\tau_{k}).
\]
Therefore, the filtering equation \ (\ref{eq:po})\ {is simply the
forward Kolmogorov equation for }$\theta.$
\end{example}

A similar effect appears also in the following example.

\begin{example}
\label{ex:constep}(Fixed observation intervals.)
Assume for simplicity that the Markov process $\theta_{t}$
is homogeneous. Also assume that $\tau_{k}=kh,$ where $h$ is a
fixed time step. Notice that
\begin{equation*}
\mathcal{G}_{t}=\mathcal{G}(k)  \ \text{for any $t\in
[\hskip-.015in[\tau_{k},\tau_{k+1}[\hskip-.015in[$}.
\end{equation*}
Denote by $P(t,x,dy)$ the transition probability kernel of the
process $\theta_{t}$, given that $\theta_{0}=x$, and let $T_{t}$
denote the associated transition operator.

In accordance with Example \ref{ex:constep0}, one can take
$\phi (\tau_{k},t)\equiv1$ and
$\Phi_{k}(dt)=\delta_{\{\tau_{k+1}\}}(t)dt.$ Thus, we get
\begin{align*}
\psi_{k}(f;t,y,\theta_{\tau_{k}})&=E\left[
f(\theta _{t})\rho_{{\tau_{k},t}}(y-X_{k})\big|\sigma\{
\theta_{\tau_{k}}\}  \vee\mathcal{G}(k)\right],
\\
\bar{\psi}_{k}(f;t,\theta_{\tau_{k}})
&=T_{t-\tau_{k} }f(\theta_{\tau_{k}})
:=\int f(y)\mathsf{P}(t-\tau_{k} ,\theta_{\tau_{k}},dy).
\end{align*}
Since $\Phi_{k}(dt)=0$ on $[\hskip-.015in[\tau_{k},\tau
_{k+1}[\hskip-.015in[$, \eqref{eq:cont} is reduced to the forward
Kolmogorov equation
\[
\frac{\partial_{t}}{\partial t}\pi_{t}(f)=\pi_{t}(\mathcal{L}f)
\]
subject to the initial condition $\pi_{\tau_{k}}(f).$ The unique
solution of this equation is given by
$\pi_{t}(f)=\pi_{\tau_{k}}(T_{t-\tau_{k}}f)$, $t<\tau_{k+1}$.
Hence,
\begin{equation*}
\pi_{\tau_{k+1}-}(f)=\pi_{\tau_{k}}(T_{h}f)
\end{equation*}
Since $\phi(\tau_{k},t)\equiv1,$ the
denominator of $\mathcal{M}_{k}$ is
equal to 1 when $t=\tau_{{k+1}}$.
This together with the formula
$\Phi(\{\tau_{k+1}\})=1$ yields
\begin{equation}
\mathcal{M}_{k}(f;t,\pi_{t})_{t=\tau_{k+1}}\Phi(\{\tau_{k+1}\})=\pi_{\tau_{k}}
\big(T_{h}f\big)-\pi_{\tau_{k+1}-}(f).
\label{eq:dopjump}
\end{equation}
Owing to (\ref{eq:dopjump}), we get
$
\mathcal{M}_{k}(f;t,\pi_{t})_{t=\tau_{k+1}}\Phi
(\{\tau_{k+1}\})=0.
$

This yields the following recursion formula:
\begin{align*}
\pi_{\tau_{k+1}}(f)&  =\frac{\pi_{\tau_{k}}(\psi _{k}(f;t,y))}{\pi_{\tau_{k}}(
\psi_{k}(1;t,y))}_{\substack{t=\tau_{k+1}\\y=X_{\tau_{k+1}}}}
\\
&  =\frac{\int_{\mathbb{R}}E\big(
f(\theta_{t-\tau_{k}})\rho_{0{,t-}\tau
_{k}}(y-z)|\theta_{0}=z\big)\pi_{\tau_{k}}(dz)}
{\int_{\mathbb{R}}E\big(
\rho_{0{,t-}\tau_{k}}(y-z)|\theta_{0}=z\big)
\pi_{\tau_{k}}\left(
dz\right)}_{\substack{t=\tau_{k+1}\\y=X_{\tau_{k+1}}}}.
\end{align*}
\end{example}

\textbf{Acknowledgement:} We are grateful to the anonymous
Associate Editor and the referee for their constructive
suggestions, especially regarding a simplified presentation of the
results. We are very much indebted to Remigijus Mikulevicius for
many important suggestions, and to Ilya Zaliapin, whose numerical
experiments helped to discover an error in a preprint version of
the paper.

\end{document}